\documentclass{amsart}
\usepackage{amsmath,amssymb}
\newtheorem{thm}{Theorem}[section]
\newtheorem{cor}[thm]{Corrolary}
\newtheorem{lem}[thm]{Lemma}

\newtheorem{propdef}[thm]{Definition-Proposition}
\newtheorem{lemdef}[thm]{Definition-Lemma}

\theoremstyle{definition}

\theoremstyle{remark}
\newtheorem{rem}[thm]{Remark}
\newtheorem{rems}[thm]{Remarks}
\numberwithin{equation}{section}

\newcommand{\co}{\colon}
\newcommand{\al}{\alpha}
\newcommand{\be}{\beta}
\newcommand{\ga}{\gamma}
\newcommand{\ali}{{\alpha^{-1}}}
\newcommand{\bei}{{\beta^{-1}}}

\newcommand{\pa}{\varphi_{\alpha}}

\newcommand{\pb}{\varphi_{\beta}}

\newcommand{\p}{$\pi$\nobreakdash-\hspace{0pt}}
\newcommand{\kt}{$\Bbbk$\nobreakdash-\hspace{0pt}}

\newcommand{\g}{\mathfrak{g}}

\newcommand{\cp}[2]{\Delta_{#1,#2}}
\newcommand{\id}{\mathrm{id}}

\newcommand{\hua}{h_{(1,\alpha)}}

\newcommand{\hdb}{h_{(2,\beta)}}

\newcommand{\htg}{h_{(3,\gamma)}}

\newcommand{\kk}{\Bbbk}
\newcommand{\cop}{\mathrm{cop}}

\newcommand{\C}{\mathbb{C}}
\newcommand{\N}{\mathbb{N}}

\newcommand{\hotimes}{\hat{\otimes}}
\newcommand{\ds}{\displaystyle}

\newcommand{\auth}{\mathrm{Aut}_\mathrm{Hopf}}

\newcommand{\R}{$R$\nobreakdash-\hspace{0pt}}
\newcommand{\h}{$h$\nobreakdash-\hspace{0pt}}
\begin{document}

\title{Graded quantum groups}

\author{Alexis Virelizier}

\address{D\'epartement des Sciences Math\'ematiques\\
Universit\'e Montpellier II \\ Case Courrier 051 \\Place Eug\`ene Bataillon\\
34095 Montpellier Cedex 5\\ France}

\email{virelizi@math.univ-montp2.fr}


\date{\today}

\begin{abstract}
Starting from a Hopf algebra endowed with an action of a group $\pi$ by Hopf automorphisms, we construct (by a
``twisted'' double method) a quasitriangular Hopf \p coalgebra. This method allows us to obtain non-trivial examples of
quasitriangular Hopf \p coalgebras for any finite group $\pi$ and for infinite groups $\pi$ such as
$\mathrm{GL}_n(\Bbbk)$. In particular, we define the graded quantum groups, which are Hopf \p coalgebras for
$\pi=\C[[h]]^l$ and generalize the Drinfeld-Jimbo quantum enveloping algebras.
\end{abstract}
\maketitle

\setcounter{tocdepth}{1} \tableofcontents

\section*{Introduction}

The aim of the present paper is to construct examples of quasitriangular Hopf group-coalgebras. These algebraic
structures were introduced by Turaev~\cite{Tur1}. Let $\pi$ be a group. The category of representations of a
quasitriangular Hopf \p coalgebra is a braided \p category. Such categories are used in~\cite{Tur1} to construct
3-dimensional homotopy quantum field theory with target an Eilenberg-Mac Lane space of type $K(\pi,1)$. Moreover, Hopf \p
coalgebras are used in \cite{mathese,Vir3} to construct Hennings-type and Kuperberg-type invariants of flat \p bundles
over link complements and over 3-manifolds.

Let $\pi$ be a group. A Hopf \p coalgebra is a family $H=\{H_\al\}_{\al \in \pi}$ of algebras (over a field $\Bbbk$)
endowed with a comultiplication $\Delta=\{\cp{\al}{\be}\co H_{\al \be} \to H_\al \otimes H_\be\}_{\al,\be \in \pi}$, a
counit $\varepsilon: H_1 \to \kk$, and an antipode $S=\{S_\al \co H_\al \to H_\ali\}_{\al \in \pi}$ which verify some
compatibility conditions. A crossing for $H$ is a family of algebra isomorphisms $\varphi=\{\varphi_\be \co H_\al \to
H_{\be \al \bei} \}_{\al,\be \in \pi}$ which preserves the comultiplication and the counit, and which yields an action of
$\pi$ in the sense that $\pb \varphi_{\be'} = \varphi_{\be \be'}$. A crossed Hopf \p coalgebra $H$ is quasitriangular
when it is endowed with an \R matrix $R=\{R_{\al,\be} \in H_\al \otimes H_\be \}_{\al,\be \in \pi}$ verifying some axioms
(involving the crossing $\varphi$) which generalize the classical ones given in~\cite{Drin2}. The case $\pi=1$ is the
standard setting of Hopf algebras.

Starting from a crossed Hopf \p coalgebra $H=\{H_\al\}_{\al \in \pi}$, Zunino \cite{Zu} constructed a double
$Z(H)=\{Z(H)_\al\}_{\al \in \pi}$ of $H$ which is a quasitriangular Hopf \p coalgebra  in which $H$ is embedded. One has
that $Z(H)_\al=H_\al \otimes \bigl (\oplus_{\be \in \pi} H_\be^* \bigr )$ as a vector space. Unfortunately, each
component $Z(H)_\al$ is infinite-dimensional (unless $H_\be=0$ for all but a finite number of $\be \in \pi$).

To obtain non-trivial examples of quasitriangular Hopf \p coalgebras with finite-dimensional components, we restrict
ourself to a less general situation: our initial data is not any crossed Hopf \p coalgebra but a Hopf algebra endowed
with an action of $\pi$ by Hopf algebra automorphisms. Remark indeed that the component $H_1$ of a Hopf \p coalgebra
$H=\{H_\al\}_{\al \in \pi}$ is a Hopf algebra and that a crossing for $H$ induces an action of $\pi$ on $H_1$ by Hopf
automorphisms.

In this paper, starting from a Hopf algebra $A$ endowed with an action $\phi \co \pi \to \auth(A)$ of a group $\pi$ by
Hopf automorphisms, we construct a quasitriangular Hopf \p coalgebra $D(A,\phi)=\{D(A,\phi_\al)\}_{\al \in \pi}$. The
algebra $D(A,\phi_\al)$ is constructed in a manner similar to the Drinfeld double (in particular $D(A,\phi_\al)=A \otimes
A^*$ as a vector space) except that its multiplication is ``twisted'' by the Hopf automorphism $\phi_\al: A \to A$. The
algebra $D(A,\id_A)$ is the usual Drinfeld double. In general, the algebras $D(A,\phi_\al)$ and $D(A,\phi_\be)$ are not
isomorphic when $\al\neq\be$.

This method allows us to define non-trivial examples of quasitriangular Hopf \p coalgebras for any finite group $\pi$ and
for infinite groups $\pi$ such as $\mathrm{GL}_n(\Bbbk)$. In particular, given a complex simple Lie algebra $\g$ of rank
$n$, we define the \emph{graded quantum groups} $\{U_h^\al(\g)\}_{\al \in \C^{*n}}$ and $\{U_h^\al(\g)\}_{\al \in
\C[[h]]^n}$ which are crossed Hopf group-coalgebras. They are obtained as quotients of $D(U_q(\mathfrak{b}_+),\phi)$ and
$D(U_h(b_+),\phi')$, where $\mathfrak{b}_+$ denotes the Borel subalgebra of $\g$, $\phi$ is an action of $\C^{*n}$ by
Hopf automorphisms of $U_q(\mathfrak{b}_+)$, and $\phi'$ is an action of $\C[[h]]^n$ by Hopf automorphisms of
$U_h(\mathfrak{b}_+)$. Furthermore, the crossed Hopf $\C[[h]]^n$-coalgebra $\{U_h^\al(\g)\}_{\al \in \C[[h]]^n}$ is
quasitriangular.

The paper is organized as follows. In Section~\ref{sect-1}, we review the basic definitions and properties of Hopf \p
coalgebras. In Section~\ref{sect-double}, we define the twisted double of a Hopf algebra $A$ endowed with an action of a
group $\pi$ by Hopf automorphisms. In Section~\ref{sect-finitegroup}, we explore the case $A=\Bbbk[G]$ where $G$ is a
finite group. In Section~\ref{sect-gnlk}, we give an example of a quasitriangular Hopf
$\mathrm{GL}_n(\Bbbk)$-coalgebra. Finally, we define the graded quantum groups in Sections~\ref{sect-uq} and~\ref{sect-uh}.\\

Throughout this paper, we let $\pi$ be a group (with neutral element~$1$) and $\Bbbk$ be a field.

\section{Hopf group-coalgebras}\label{sect-1}
In this section, we review some definitions and properties concerning Hopf group-coalgebras. For a detailed treatment of
the theory of Hopf group-coalgebras, we refer to \cite{Vir2}.

\subsection{Hopf \protect\p coalgebras}\label{hopfpicoal}
A \emph{Hopf \p coalgebra} (over $\Bbbk$) is a family $H=\{H_\al\}_{\al \in \pi}$ of \kt algebras endowed with a family
$\Delta=\{\cp{\al}{\be}\co H_{\al \be} \to H_\al \otimes H_\be \}_{\al,\be \in \pi}$ of algebra homomorphisms (the
\emph{comultiplication}) and an algebra homomorphism $\varepsilon\co H_1 \to \Bbbk$ (the \emph{counit}) such that, for
all $\al,\be,\ga \in \pi$,
\begin{align}
& (\cp{\al}{\be}\otimes \id_{H_\ga}) \cp{\al \be}{\ga}=(\id_{H_\al} \otimes \cp{\be}{\ga}) \cp{\al}{\be \ga},
  \label{coass} \\
& (\id_{H_\al} \otimes \varepsilon) \cp{\al}{1}=\id_{H_\al}=(\varepsilon \otimes \id_{H_\al}) \cp{1}{\al} \label{counit},
\end{align}
and with a family $S=\{S_\al \co  H_\al \to H_\ali \}_{\al \in \pi}$ of $\Bbbk$-linear maps (the \emph{antipode}) which
verifies that, for all $\al \in \pi$,
\begin{equation}\label{antipode}
m_\al (S_\ali \otimes \id_{H_\al}) \cp{\ali}{\al} = \varepsilon \,1_\al =  m_\al (\id_{H_\al} \otimes S_\ali)
\cp{\al}{\ali},
\end{equation}
where $m_\al\co H_\al \otimes H_\al \to H_\al$ and $1_\al \in H_\al$ denote respectively the multiplication and unit
element of $H_\al$.

When $\pi=1$, one recovers the usual notion of a Hopf algebra. In particular $(H_1,m_1,1_1,\cp{1}{1},\varepsilon,S_1)$ is
a Hopf algebra.

Remark that the notion of a Hopf \p coalgebra is not self-dual and that if $H=\{H_\al\}_{\al \in \pi}$ is a Hopf \p
coalgebra, then $\{ \al \in \pi \, | \, H_\al \neq 0 \}$ is a subgroup of~$\pi$.

A Hopf \p coalgebra $H=\{H_\al\}_{\al \in \pi}$ is said to be of \emph{finite type} if, for all $\al \in \pi$, $H_\al$ is
finite-dimensional (over $\Bbbk$). Note that it does not mean that $\oplus_{\al \in \pi} H_\al$ is finite-dimensional
(unless $H_\al=0$ for all but a finite number of $\al \in \pi$).

The antipode of a Hopf \p coalgebra $H=\{H_\al\}_{\al \in \pi}$ is anti-multiplicative: each $S_\al\co  H_\al \to H_\ali$
is an anti-homomorphism of algebras, and anti-comultiplicative: $\varepsilon S_1=\varepsilon$ and $\cp{\bei}{\ali} S_{\al
\be}=\sigma_{H_\ali,H_\bei}(S_\al \otimes S_\be) \cp{\al}{\be}$ for any $\al,\be \in \pi$, see \cite[Lemma~1.1]{Vir2}.

The antipode $S=\{ S_\al \}_{\al \in \pi}$ of $H=\{H_\al\}_{\al \in \pi}$ is said to be \emph{bijective} if each $S_\al$
is bijective. As for Hopf algebras, the antipode of a finite type Hopf \p coalgebra is always bijective (see
\cite[Corollary~3.7(a)]{Vir2}).

We extend the Sweedler notation for the comultiplication of a Hopf \p coalgebra $H=\{H_\al\}_{\al \in \pi}$ in the
following way: for any $\al,\be \in \pi$ and $h \in H_{\al \be}$, we write $
  \cp{\al}{\be}(h)=\sum_{(h)} \hua \otimes \hdb \in H_\al \otimes H_\be,
$ or shortly, if we leave the summation implicit,
  $
  \cp{\al}{\be}(h)= \hua \otimes \hdb
  $. The coassociativity of $\Delta$ gives that, for any $\al,\be, \ga \in \pi$ and $h \in H_{\al \be \ga}$,
  \begin{equation*}
  h_{(1,\al \be)(1,\al)} \otimes h_{(1,\al \be)(2,\be)} \otimes
  h_{(2,\ga)} = h_{(1,\al)} \otimes h_{(2,\be \ga)(1,\be)} \otimes h_{(2,\be \ga)(2,\ga)}.
  \end{equation*}
This element of $H_\al \otimes H_\be \otimes H_\ga$ is written as $\hua \otimes \hdb \otimes \htg$. By iterating the
procedure, we define inductively $h_{(1,\al_1)} \otimes \dotsb \otimes h_{(n,\al_n)}$ for any $h \in H_{\al_1 \dotsm
\al_n}$.

\subsection{Crossed Hopf \protect\p coalgebras}\label{deficro}
A Hopf \p coalgebra $H=\{H_\al\}_{\al \in \pi}$ is said to be \emph{crossed} if it is endowed with a family
$\varphi=\{\varphi_\be \co  H_\al \to H_{\be \al \bei} \}_{\al,\be \in \pi}$ of algebra isomorphisms (the
\emph{crossing}) such that, for all $\al,\be,\ga \in \pi$,
\begin{align}
& (\varphi_\be \otimes \varphi_\be) \cp{\al}{\ga} =\cp{\be \al \bei}{\be \ga \bei} \varphi_\be, \label{Crossdef1}\\
& \varepsilon \varphi_\be=\varepsilon, \label{Crossdef2}\\
& \varphi_\al \varphi_\be=\varphi_{\al \be}. \label{Crossdef3}
\end{align}

It is easy to check that $\varphi_1|_{H_\al}=\id_{H_\al}$ and $\varphi_\be S_\al = S_{\be \al \bei} \varphi_\be$ for all
$\al,\be \in \pi$.

\subsection{Quasitriangular Hopf \protect\p coalgebras}\label{quasitrig}
A crossed Hopf \p coalgebra $H=\{H_\al\}_{\al \in \pi}$ is said to be \emph{quasitriangular} if it is endowed with  a
family $R=\{R_{\al,\be} \in H_\al \otimes H_\be\}_{\al,\be \in \pi}$ of invertible elements (the \emph{\R matrix}) such
that, for all $\al,\be, \ga \in \pi$ and $x \in H_{\al\be}$,
\begin{align}
 & R_{\al,\be} \cdot \cp{\al}{\be}(x)= \sigma_{\be,\al}
             (\varphi_\ali\otimes \id_{H_\al}) \cp{\al \be \ali}{\al}(x) \cdot R_{\al,\be}, \label{QT1}\\
 & (\id_{H_\al} \otimes \cp{\be}{\ga})(R_{\al,\be \ga})=(R_{\al,\ga})_{1 \be 3}
     \cdot (R_{\al,\be})_{12 \ga}, \label{QT2} \\
 & (\cp{\al}{\be} \otimes \id_{H_\ga})(R_{\al \be, \ga})=[(\id_{H_\al} \otimes
     \varphi_{\bei})(R_{\al, \be \ga \bei})]_{1 \be 3} \cdot (R_{\be,\ga})_{\al 23},\label{QT3} \\
 & (\varphi_\be \otimes \varphi_\be)(R_{\al,\ga})=R_{\be \al \bei, \be \ga \bei}, \label{QT4}
\end{align}
where $\sigma_{\be,\al}$ denotes the flip map $H_\be \otimes H_\al \to H_\al \otimes H_\be$ and, for \kt spaces $P,Q$ and
$r=\sum_j p_j \otimes q_j \in P \otimes Q$, we set $r_{12\ga}=r \otimes 1_\ga \in P \otimes Q \otimes H_\ga$, $r_{\al
23}=1_\al \otimes r \in H_\al \otimes P \otimes Q$, and $r_{1 \be 3}=\sum_j p_j \otimes 1_\be \otimes q_j \in P \otimes
H_\be \otimes Q$.

Note that $R_{1,1}$ is a (classical) \R matrix for the Hopf algebra $H_1$.

When $\pi$ is abelian and $\varphi$ is \emph{trivial} (that is, $\varphi_\be |_{H_\al}=\id_{H_\al}$ for all $\al,\be \in
\pi$), one recovers the definition of a quasitriangular \p colored Hopf algebra given by Ohtsuki in~\cite{Oh}.

The \R matrix always verifies (see \cite[Lemma~6.4]{Vir2}) that, for any $\al, \be, \ga \in \pi$,
\begin{align}
  & (\varepsilon \otimes \id_{H_\al} )(R_{1,\al}) = 1_\al =(\id_{H_\al} \otimes
        \varepsilon) (R_{\al,1}), \\
  & (S_\ali \varphi_\al \otimes \id_{H_\be} )(R_{\ali,\be})= R_{\al,\be}^{-1} \text{\quad and \quad}
    (\id_{H_\al} \otimes S_\be)(R_{\al,\be}^{-1})= R_{\al,\bei}, \\
  & (S_\al \otimes S_\be)(R_{\al,\be})= (\varphi_\al \otimes
        \id_{H_{\be^{-1}}})(R_{\al^{-1},\be^{-1}}),
\end{align}
and provides a solution of the \p colored Yang-Baxter equation:
\begin{equation}
  \begin{split} (R_{\be,\ga})_{\al 23} & \cdot (R_{\al,\ga})_{1 \be 3}  \cdot (R_{\al,\be})_{12 \ga}\\
      & =(R_{\al,\be})_{1 2 \ga} \cdot [(\id_{H_\al} \otimes \varphi_\bei)
        (R_{\al,\be \ga \bei})]_{1 \be 3}\cdot (R_{\be,\ga})_{\al 23}.
      \end{split}
\end{equation}

\subsection{Ribbon Hopf \protect\p coalgebras}\label{ribpico}
A quasitriangular Hopf \p coalgebra $H=\{H_\al\}_{\al \in \pi}$ is said to be \emph{ribbon} if it is endowed with a
family $\theta=\{\theta_\al \in H_\al \}_{\al \in \pi}$ of invertible elements (the \emph{twist}) such that, for any
$\al,\be \in \pi$,
\begin{align}
  & \varphi_\al(x)=\theta_\ali x \theta_\al \text{\quad for all $ x \in H_\al$},\label{twist1} \\
  & S_\al(\theta_\al)=\theta_\ali, \label{twist2} \\
  & \pb(\theta_\al)=\theta_{\be \al \bei}, \label{twist3} \\
  & \cp{\al}{\be}(\theta_{\al \be})=(\theta_\al \otimes \theta_\be) \cdot
            \sigma_{\be,\al}( (\varphi_\ali \otimes \id_{H_\al} )(R_{\al \be \ali ,\al}))
            \cdot R_{\al,\be}.\label{twist4}
\end{align}

Note that $\theta_1$ is a (classical) twist of the quasitriangular Hopf algebra $H_1$.

\subsection{Hopf $\pi$-coideals}
Let $H=\{H_\al\}_{\al \in \pi}$ be a Hopf \p coalgebra. A \emph{Hopf \p coideal} of $H$ is a family $I=\{I_\al\}_{\al \in
\pi}$, where each $I_\al$ is an ideal of $H_\al$, such that, for any $\al,\be \in \pi$,
\begin{align}
 & \cp{\al}{\be}(I_{\al \be}) \subset I_\al \otimes H_\be + H_\al \otimes I_\be,\\
 & \varepsilon(I_1)=0,\\
 & S_\al(I_\al) \subset I_\ali.
\end{align}
The quotient $\Bar{H}=\{ \Bar{H}_\al=H_\al / I_\al \}_{\al \in \pi}$, endowed with the induced structure maps, is then a
Hopf \p coalgebra. If $H$ is furthermore crossed, with a crossing $\varphi$ such that, for any $\al,\be \in \pi$,
\begin{equation}
\pb(I_\al) \subset I_{\al\be \ali},
\end{equation}
then so is $\Bar{H}$ (for the induced crossing).

\section{Twisted double of Hopf algebras}\label{sect-double}

In this section, we give a method (the \emph{twisted double}) of defining a quasitriangular Hopf \p coalgebra from a Hopf
algebra endowed with an action of a group $\pi$ by Hopf automorphisms.

\subsection{Hopf pairings}\label{hopfpar}
Recall that a \emph{Hopf pairing} between two Hopf algebras $A$ and $B$ (over $\Bbbk$) is a bilinear pairing $\sigma \co
A \times B \to \Bbbk$ such that, for all $a, a' \in A$ and $b,b' \in B$,
\begin{align}
 & \sigma ( a, b b' ) = \sigma ( a_{(1)}, b)\, \sigma ( a_{(2)}, b') , \\
 & \sigma ( a a', b) = \sigma ( a, b_{(2)}) \,\sigma ( a', b_{(1)}) , \\
 & \sigma ( a, 1 )=\varepsilon(a) \text{\quad and \quad }\sigma ( 1, b )=\varepsilon(b).
\end{align}
Note that such a pairing always verifies that, for any $a\in A$ and $b\in B$,
\begin{align}
 & \sigma ( S(a), S(b) )=\sigma ( a, b ).
\end{align}
(Since both $\sigma$ and $\sigma \circ (S\times S)$ are the inverse of $\sigma \circ (\id \times S)$ in the algebra
$\mathrm{Hom}_\kk(A \times B, \kk)$ endowed with the convolution product).\\

Let $\sigma \co A \times B \to \Bbbk$ be a Hopf pairing. Its annihilator ideals are $ I_A=\{a \in A \, | \, \sigma(a,b)=0
\text{ for all } b \in B \} $ and $ I_B=\{b \in B \, | \, \sigma(a,b)=0 \text{ for all } a \in A \} $. It is easy to
check that $I_A$ and $I_B$ are Hopf ideals of $A$ and $B$, respectively. Recall that $\sigma$ is said to be
\emph{non-degenerate} if $I_A$ and $I_B$ are both reduced to $0$. A degenerate Hopf pairing $\sigma \co A \times B \to
\Bbbk$ induces
(by passing to the quotients) a Hopf pairing $\Bar{\sigma} \co A/I_A \times B/I_B \to \Bbbk$ which is non-degenerate.\\

Most of Hopf algebras we shall consider in the sequel will be defined by generators and relations. The following provides
us with a method of constructing Hopf pairings, see \cite{Dae,KRT}.

Let $\Tilde{A}$ (resp.\@ $\Tilde{B}$) be a free algebra generated by elements $a_1, \dots,a_p$ (resp.\@ $b_1,\dots ,
b_q$) over $\Bbbk$. Suppose that $\Tilde{A}$ and $\Tilde{B}$ have Hopf algebra structures such that each $\Delta(a_i)$
for $1\leq i \leq p$ (resp.\@ $\Delta(b_j)$ for $1\leq i \leq q$) is a linear combination of tensors $a_r \otimes a_s$
(resp.\@ $b_r \otimes b_s$). Given $pq$ scalars $\lambda_{i,j} \in \Bbbk$ with $1 \leq i \leq p$ and $1 \leq j \leq q$,
there is a unique Hopf pairing $\sigma: \Tilde{A} \times \Tilde{B} \to \Bbbk$ such that $\sigma(a_i,b_j)=\lambda_{i,j}$.

Suppose now that $A$ (resp.\@ $B$) is the algebra obtained as the quotient of $\Tilde{A}$ (resp.\@ $\Tilde{B}$) by the
ideal generated by elements $r_1, \dots,r_m \in \Tilde{A}$ (resp.\@ $s_1,\dots,s_n \in \Tilde{B}$). Suppose also that the
Hopf algebra structure in $\Tilde{A}$ (resp.\@ $\Tilde{B}$) induces a Hopf algebra structure in $A$ (resp.\@ $B$). Then a
Hopf pairing $\sigma: \Tilde{A} \times \Tilde{B} \to \Bbbk$ induces a Hopf pairing $A \times B \to \Bbbk$ if and only if
$\sigma(r_i,b_j)=0$ for all $1\leq i \leq m$ and $1\leq j \leq q$, and $\sigma(a_i,s_j)=0$ for all $1\leq i \leq p$ and
$1\leq j \leq n$.

\subsection{The twisted double construction}

\begin{lemdef}\label{algdouble}
Let $\sigma \co A \times B \to \Bbbk$ be Hopf pairing between two Hopf algebras $A$ and $B$. Let $\phi \co A \to A$ be a
Hopf algebra endomorphism of $A$. Set $D(A,B;\sigma,\phi)=A \otimes B$ as a \kt space. Then $D(A,B;\sigma,\phi)$ has a
structure of an associative and unitary algebra given, for any $a,a' \in A$ and $b,b' \in B$, by
\begin{align}
 & (a \otimes b) \cdot (a'\otimes b')=\sigma ( \phi(a'_{(1)}),S(b_{(1)}) ) \, \sigma ( a'_{(3)},b_{(3)} ) \; a a'_{(2)}
   \otimes b_{(2)} b', \label{crossrel1}\\
 & 1_{D(A,B;\sigma,\phi)}=1_A \otimes 1_B.
\end{align}
Moreover, the linear embeddings $A\hookrightarrow D(A,B;\sigma,\phi)$ and $B\hookrightarrow D(A,B;\sigma,\phi)$ defined
by $a \mapsto a \otimes 1_B$ and $b \mapsto 1_A \otimes b$, respectively, are algebra morphisms.
\end{lemdef}
\begin{rems}  {(a)} Note that $D(A,B;\sigma,\id_A)$ is the underlying algebra of the usual quantum double of $A$ and $B$
(obtained by
using the Hopf pairing $\sigma$).\\
\indent {(b)} If $\phi$ and $\phi'$ are different Hopf algebra endomorphisms of $A$, then the algebras
$D(A,B;\sigma,\phi)$ and $D(A,B;\sigma,\phi')$ are not in general isomorphic, see Remark~\ref{remnoteq}.
\end{rems}
\begin{proof}
Let $a,a',a'' \in A$ and $b,b',b'' \in B$. Using the fact that $\sigma$ is a Hopf pairing and $\phi$ is a Hopf algebra
endomorphism, we have that
\begin{eqnarray*}
\lefteqn{((a\otimes b) \cdot  (a' \otimes b')) \cdot (a'' \otimes b'') } &&\\
 & = & \sigma(\phi(a'_{(1)}), S(b_{(1)})) \, \sigma(a'_{(3)},b_{(5)}) \, \sigma(\phi(a''_{(1)}), S(b_{(2)}b'_{(1)})) \\
 &   & \phantom{xxxxxxxxxxxxxx}
       \sigma(a''_{(3)}, b_{(4)}b'_{(3)}) \; aa'_{(2)}a''_{(3)} \otimes b_{(3)}b'_{(2)}b''\\
 & = & \sigma(\phi(a'_{(1)}), S(b_{(1)})) \, \sigma(a'_{(3)},b_{(5)}) \, \sigma(\phi(a''_{(1)}), S(b'_{(1)}))
       \, \sigma(\phi(a''_{(2)}), S(b_{(2)})) \\
 &   & \phantom{xxxxxxxxxxxxxx}
       \sigma(a''_{(4)}, b_{(4)}) \, \sigma(a''_{(5)},b'_{(3)}) \; aa'_{(2)}a''_{(3)} \otimes b_{(3)}b'_{(2)}b'',
\end{eqnarray*}
and
\begin{eqnarray*}
\lefteqn{  (a\otimes b) \cdot ((a' \otimes b')   \cdot (a'' \otimes b''))} && \\
& = & \sigma(\phi(a''_{(1)}), S(b'_{(1)})) \, \sigma(a''_{(5)},b'_{(3)}) \, \sigma(\phi(a'_{(1)}a''_{(2)}), S(b_{(1)})) \\
 &   & \phantom{xxxxxxxxxxxxxx}
       \sigma(a'_{(3)} a''_{(4)}, b_{(3)}) \; aa'_{(2)}a''_{(3)} \otimes b_{(2)}b'_{(2)}b''\\
 & = & \sigma(\phi(a''_{(1)}), S(b'_{(1)})) \, \sigma(a''_{(5)},b'_{(3)}) \, \sigma(\phi(a'_{(1)}), S(b_{(1)}))
       \, \sigma(\phi(a''_{(2)}), S(b_{(2)})) \\
 &   & \phantom{xxxxxxxxxxxxxx}
       \sigma(a'_{(3)}, b_{(5)}) \, \sigma(a''_{(4)},b_{(4)}) \; aa'_{(2)}a''_{(3)} \otimes b_{(3)}b'_{(2)}b''.
\end{eqnarray*}
Hence the product is associative. Finally, $1_A \otimes 1_B$ is the unit element since
\begin{eqnarray*}
(a\otimes b) \cdot (1 \otimes 1)& = &  \sigma ( \phi(1),S(b_{(1)}) ) \, \sigma ( 1,b_{(3)} ) \; a
\otimes b_{(2)}\\
 & = & \varepsilon(S(b_{(1)}) ) \, \varepsilon(b_{(3)} ) \; a \otimes b_{(2)} \;\, = \;\, a \otimes b,
\end{eqnarray*}
and
\begin{eqnarray*}
(1 \otimes 1) \cdot (a\otimes b) & = &  \sigma ( \phi(a_{(1)}),S(1) ) \, \sigma ( a_{(3)},1 ) \; a_{(2)} \otimes b\\
 & = & \varepsilon(\phi(a_{(1)})) \, \varepsilon(a_{(3)} ) \; a_{(2)} \otimes b \;\, = \;\, a \otimes b.
\end{eqnarray*}
Finally, for any $a,a' \in A$ and $b,b' \in B$, we have that
\begin{eqnarray*}
 (a \otimes 1) \cdot (a'\otimes 1) & = & \sigma ( \phi(a'_{(1)}),S(1) ) \, \sigma ( a'_{(3)},1 ) \; a a'_{(2)}
   \otimes 1 \\
 & = & \varepsilon(\phi(a'_{(1)})) \, \varepsilon(a'_{(3)}) \; a a'_{(2)} \otimes 1 \\
 & = &  a a' \otimes 1
\end{eqnarray*}
and
\begin{eqnarray*}
 (1 \otimes b) \cdot (1\otimes b') & = & \sigma ( \phi(1),S(b_{(1)}) ) \, \sigma (1, b_{(3)}) \; 1 \otimes b_{(2)}b'\\
 & = & \varepsilon(S(b_{(1)})) \, \varepsilon(b_{(3)}) \; 1 \otimes b_{(2)}b' \\
 & = &  1 \otimes b b'.
\end{eqnarray*}
Therefore $A\hookrightarrow D(A,B;\sigma,\phi)$ and $B\hookrightarrow D(A,B;\sigma,\phi)$ are algebra morphisms.
\end{proof}

In the sequel, the group of Hopf automorphisms of a Hopf algebra $A$ will be denoted by $\auth(A)$.

\begin{thm}\label{Hopfdouble}
Let $\sigma \co A \times B \to \Bbbk$ be Hopf pairing between two Hopf algebras $A$ and $B$, and $\phi \co \pi \to \auth
(A)$ be group homomorphism (that is, an action of $\pi$ on $A$ by Hopf automorphisms). Then the family of algebras
$D(A,B;\sigma, \phi)=\{ D(A,B;\sigma, \phi_\al)\}_{\al \in \pi}$ (see Definition~\ref{algdouble}) has a structure of a
Hopf \p coalgebra given, for any $a\in A$, $b \in B$, and $\al,\be \in \pi$, by
\begin{align}
& \Delta_{\alpha,\beta}(a \otimes b)=(\phi_\beta(a_{(1)}) \otimes b_{(1)} ) \otimes (a_{(2)} \otimes b_{(2)}), \\
& \varepsilon(a \otimes b)=\varepsilon_A(a) \, \varepsilon_B(b), \\
& S_\alpha (a \otimes b) =\sigma(\pa(a_{(1)}),b_{(1)}) \, \sigma(a_{(3)},S(b_{(3)})) \; \pa S(a_{(2)}) \otimes
S(b_{(2)}).
\end{align}
\end{thm}
\begin{proof}
The coassociativity \eqref{coass} follows directly from the coassociativity of the coproducts of $A$ and $B$ and the fact
that $\phi_{\be\ga}=\phi_\be \circ \phi_\ga$. Axiom \eqref{counit} is a direct consequence of $\varepsilon_A \circ
\phi_\al =\varepsilon_A$. Since $\phi_1=\id_A$ and $D(A,B;\sigma,\id_A)$ is underlying algebra of the usual quantum
double of $A$ and $B$, the counit $\varepsilon$ is multiplicative. Let us verify that $\Delta_{\alpha,\beta}$ is
multiplicative. Let $a,a' \in A$ and $b,b' \in B$. On one hand we have:
\begin{eqnarray*}
\lefteqn{\cp{\al}{\be}((a \otimes b) \cdot (a' \otimes b'))} && \\
  & = & \sigma(\phi_{\al\be}(a'_{(1)}),S(b_{(1)})) \, \sigma(a'_{(3)}, b_{(3)}) \; \cp{\al}{\be}(a a'_{(2)} \otimes b_{(2)} b') \\
  & = & \sigma(\phi_{\al\be}(a'_{(1)}),S(b_{(1)})) \, \sigma(a'_{(4)}, b_{(4)}) \; \phi_\be(a_{(1)}a'_{(2)})
        \otimes b_{(2)} b'_{(1)} \otimes a_{(2)}a'_{(3)} \otimes b_{(3)}b'_{(2)}.
\end{eqnarray*}
One the other hand,
\begin{eqnarray*}
\lefteqn{\cp{\al}{\be}(a \otimes b) \cdot \cp{\al}{\be}(a' \otimes b')} &&\\
  & = & (\phi_\be (a_{(1)}) \otimes b_{(1)} \otimes a_{(2)} \otimes b_{(2)}) \cdot (\phi_\be(a'_{(1)}) \otimes b'_{(1)}
        \otimes a'_{(2)} \otimes b'_{(2)}) \\
  & = & \sigma(\phi_\al \phi_\be(a'_{(1)}), S(b_{(1)})) \, \sigma(\phi_\be(a'_{(3)}),b_{(3)}) \, \sigma
        (\phi_\be(a'_{(4)}),S(b_{(4)})) \, \sigma(a'_{(6)},b_{(6)}) \\
  &   &  \phantom{xxxxxxxxxxxxxx} \phi_\be(a_{(1)})\phi_\be(a'_{(2)}) \otimes b_{(2)} b'_{(1)}
        \otimes a_{(2)} a'_{(5)} \otimes b_{(5)} b'_{(2)} \\
  & = & \sigma(\phi_{\al\be}(a'_{(1)}), S(b_{(1)})) \, \sigma(\phi_\be(a'_{(3)}),b_{(3)}S(b_{(4)}))
        \, \sigma(a'_{(5)},b_{(6)}) \\
  &   &  \phantom{xxxxxxxxxxxxxx} \phi_\be(a_{(1)}a'_{(2)}) \otimes b_{(2)} b'_{(1)}
        \otimes a_{(2)} a'_{(4)} \otimes b_{(5)} b'_{(2)} \\
  & = & \sigma(\phi_{\al \be}(a'_{(1)}), S(b_{(1)})) \, \sigma(a'_{(4)},b_{(4)})\phi_\be(a_{(1)} a'_{(2)})
        \otimes b_{(2)} b'_{(1)}
        \otimes a_{(2)} a'_{(3)} \otimes b_{(3)} b'_{(2)}.
\end{eqnarray*}
Let us verify the first equality of \eqref{antipode}. Let $a \in A$, $b \in B$, and $\al \in \pi$. Denote the
multiplication in $D(A,B;\sigma,\phi_\al)$ by $m_\al$. We have
\begin{eqnarray*}
\lefteqn{m_\al(S_\ali \otimes \id_{D(A,B;\sigma,\phi_\al)})\cp{\ali}{\al}(a \otimes b)} &&\\
  & = & \sigma(a_{(1)},b_{(1)}) \,\sigma (\phi_\al(a_{(3)}),S(b_{(5)})) \, \sigma(\phi_\al(a_{(4)}), S^2(b_{(4)})) \,
        \sigma(a_{(6)},S(b_{(2)})) \\
  &   & \phantom{xxxxxxxxxxxxxx} S(a_{(2)})a_{(5)} \otimes S(b_{(3)})b_{(6)} \\
  & = & \sigma(a_{(1)},b_{(1)}) \,\sigma (\phi_\al(a_{(3)}),S(b_{(5)}) S^2(b_{(4)})) \,
        \sigma(a_{(5)},S(b_{(2)})) \\
  &   & \phantom{xxxxxxxxxxxxxx} S(a_{(2)})a_{(4)} \otimes S(b_{(3)})b_{(6)} \\
  & = & \sigma(a_{(1)},b_{(1)}) \, \sigma(a_{(4)},S(b_{(2)})) \; S(a_{(2)})a_{(3)} \otimes S(b_{(3)})b_{(4)} \\
  & = & \sigma(a_{(1)},b_{(1)}) \, \sigma(a_{(2)},S(b_{(2)})) \; 1 \otimes 1 \\
  & = & \sigma(a,b_{(1)}S(b_{(2)})) \; 1 \otimes 1  \;\, = \;\, \varepsilon(a) \, \varepsilon(b) \; 1 \otimes 1
\end{eqnarray*}
The second equality of \eqref{antipode} can be verified similarly.
\end{proof}

Let $\sigma \co A \times B \to \Bbbk$ be a Hopf pairing between two Hopf algebras $A$ and $B$, and $\phi \co \pi \to
\auth (A)$ be an action of $\pi$ on $A$ by Hopf automorphisms. An action $\psi \co \pi \to \auth (B)$ of $\pi$ on $B$ by
Hopf automorphisms is said to be \emph{$(\sigma,\phi)$-compatible} if, for all $a \in A$, $b \in B$ and $\beta \in \pi$,
\begin{equation}
\sigma ( \phi_\beta(a), \psi_{\beta}(b) )=\sigma ( a, b ).
\end{equation}

\begin{lem}\label{crosseddegen}
Let $\sigma \co A \times B \to \Bbbk$ be a Hopf pairing between two Hopf algebras $A$ and $B$, and $\phi \co \pi \to
\auth (A)$, $\psi \co \pi \to \auth (B)$ be two actions of $\pi$ by Hopf automorphisms on $A$ and $B$, respectively.
Suppose that $\psi$ is $(\sigma,\phi)$-compatible. Then the Hopf \p coalgebra $D(A,B;\sigma, \phi)=\{ D(A,B;\sigma,
\phi_\al)\}_{\al \in \pi}$ (see Theorem~\ref{Hopfdouble}) admits a crossing $\varphi$ given, for any $a \in A$, $b\in B$
and $\be \in \pi$, by
\begin{equation}
\varphi_\beta(a\otimes b)=\phi_\beta (a) \otimes \psi_\beta(b).
\end{equation}
\end{lem}
\begin{proof}
Let $\al,\be \in \pi$. We have that $\pb(1_A \otimes 1_B)=\phi_\be(1_A) \otimes \psi_\be(1_B)=1_A \otimes 1_B$ and, for
any $a,a' \in A$ and $b,b' \in B$,
\begin{eqnarray*}
\lefteqn{\varphi_\be(a \otimes b)\cdot \varphi_\be(a'\otimes b')} &&\\
 & = & \sigma ( \phi_{\be\al\bei}(\phi_\be(a')_{(1)}),S(\psi_\be(b)_{(1)}) ) \, \sigma ( \phi_\be(a')_{(3)},\psi_\be(b)_{(3)} )
       \; \phi_\be(a) \phi_\be(a')_{(2)} \otimes \psi_\be (b)_{(2)}\psi_\be(b') \\
 & = & \sigma ( \phi_\be \phi_\al(a'_{(1)})),\psi_\be S(b_{(1)}) ) \, \sigma ( \phi_\be(a'_{(3)}),\psi_\be(b_{(3)}) )
       \; \phi_\be(a) \phi_\be(a'_{(2)}) \otimes \psi_\be (b_{(2)})\psi_\be(b') \\
 & = & \sigma (\phi_\al(a'_{(1)})),S(b_{(1)}) ) \, \sigma ( a'_{(3)},b_{(3)} )
       \; \phi_\be(aa'_{(2)}) \otimes \psi_\be (b_{(2)}b') \\
 & = & \varphi_\be((a \otimes b)\cdot (a'\otimes b')).
\end{eqnarray*}
Moreover $\phi_\be$ and $\psi_\be$ are bijective and so is $\varphi_\be$. Therefore $\varphi_\be : D(A,B;\sigma,
\phi_\al) \to D(A,B;\sigma, \phi_{\be \al \bei})$ is an algebra isomorphism.

Let $a \in A$, $b \in B$ and $\al,\be, \ga \in \pi$. We have that
\begin{eqnarray*}
\cp{\be \al \bei}{\be \ga \bei} (\varphi_\be(a \otimes b))
 & = & \phi_{\be \ga \bei}(\phi_\be(a)_{(1)}) \otimes \psi_\be(b)_{(1)} \otimes \phi_\be(a)_{(2)} \otimes
       \psi_\be(b)_{(2)}\\
 & = & \phi_{\be \ga \bei}\phi_\be(a_{(1)}) \otimes \psi_\be(b_{(1)}) \otimes \phi_\be(a_{(2)}) \otimes
       \psi_\be(b_{(2)})\\
 & = & \phi_\be \phi_\ga(a_{(1)}) \otimes \psi_\be(b_{(1)}) \otimes \phi_\be(a_{(2)}) \otimes
       \psi_\be(b_{(2)})\\
 & = & (\varphi_\be \otimes \varphi_\be) \cp{\al}{\ga}(a \otimes b),
\end{eqnarray*}
\begin{equation*}
\varepsilon \varphi_\be(a \otimes b)=\varepsilon(\phi_\be(a)) \, \varepsilon(\psi_\be(b)) =\varepsilon(a) \,
\varepsilon(b)=\varepsilon(a \otimes b),
\end{equation*}
and
\begin{equation*}
\varphi_\al \varphi_\be(a \otimes b)
 =\phi_\al \phi_\be(a)\otimes \psi_\al \psi_\be(b) \\
 = \phi_{\al \be}(a) \otimes \psi_{\al\be}(b)\\
 = \varphi_{\al \be}(a \otimes b).
\end{equation*}
Therefore $\varphi$ satisfies Axioms \eqref{Crossdef1}, \eqref{Crossdef2} and \eqref{Crossdef3}.
\end{proof}

\begin{cor}\label{crosednondegen}
Let $\sigma\co A \times B \to \Bbbk$ be a non-degenerate Hopf pairing and $\phi \co \pi \to \auth (A)$ be an action of
$\pi$ on $A$ by Hopf automorphisms. Then there exists a unique action $\phi^* \co \pi \to \auth (B)$ which is
$(\sigma,\phi)$-compatible. It is characterized, for any $a \in A$, $b \in B$ and $\be \in \pi$, by
\begin{equation}\label{crossednondegen}
\sigma(a,\phi_\be^*(b))=\sigma(\phi_\bei(a),b).
\end{equation}
Consequently the Hopf \p coalgebra $D(A,B;\sigma, \phi)=\{ D(A,B;\sigma, \phi_\al)\}_{\al \in \pi}$ (see
Theorem~\ref{Hopfdouble}) is crossed with crossing defined by $\varphi_\be=\phi_\be \otimes \phi^*_\be$ for any $\be \in
\pi$.
\end{cor}
\begin{proof}
Let $\be \in \pi$. Since $\sigma$ is non-degenerate, \eqref{crossednondegen} does define a linear map $\phi_\be^* \co B
\to B$. Since $\sigma$ is a Hopf pairing and $\phi_\bei$ is a Hopf algebra isomorphism of $A$, the map $\phi_\be^*$ is a
Hopf algebra isomorphism of $B$. Moreover $\phi^*$ is an action since $\phi_1^*=\id_B$ (because $\phi_1=\id_A$) and, for
any $a\in A$, $b \in B$ and $\al,\be \in \pi$,
\begin{multline*}
\sigma(a,\phi_{\al \be}^*(b))= \sigma(\phi_{\bei \ali}(a),b) =\sigma(\phi_\bei\phi_\ali(a),b) \\
=\sigma(\phi_\ali(a),\phi_\be^*(b))=\sigma(a,\phi_\al^*\phi_\be^*(b)).
\end{multline*}
Finally \eqref{crossednondegen} says exactly that $\phi^*$ is $(\sigma,\phi)$-compatible.
\end{proof}

\begin{thm}\label{quasidouble}
Let $\sigma \co A \times B \to \Bbbk$ be a Hopf pairing between two Hopf algebras $A$ and $B$, and $\phi \co \pi \to
\auth (A)$ be an action of $\pi$ on $A$ by Hopf automorphisms. Suppose that $\sigma$ is non-degenerate and that $A$ (and
so $B$) is finite dimensional. Then the crossed Hopf \p coalgebra $D(A,B;\sigma, \phi)=\{ D(A,B;\sigma, \phi_\al)\}_{\al
\in \pi}$ (see Corollary~\ref{crosednondegen}) is quasitriangular with $R$-matrix given, for all $\al,\be \in \pi$, by
\begin{equation}
R_{\al,\be}=\sum_i (e_i \otimes 1_B ) \otimes (1_A \otimes f_i),
\end{equation}
where $(e_i)_i$ and $(f_i)_i$ are basis of $A$ and $B$, respectively, such that $\sigma(e_i,f_j)=\delta_{i,j}$.
\end{thm}
\begin{rems}
{(a)} The element $\sum_i (e_i \otimes 1_B ) \otimes (1_A \otimes f_i) \in A \otimes B \otimes A \otimes B$ is canonical,
i.e., independent of the choices of the basis $(e_i)_i$ of $A$ and $(f_i)_i$ of $B$ such that
$\sigma(e_i,f_j)=\delta_{i,j}$.\\
\indent {(b)} The hypothesis $A$ is finite dimensional is to ensure that the sum $\sum_i (e_i \otimes 1_B ) \otimes (1_A
\otimes f_i)$ lies in $A \otimes B \otimes A \otimes B$. More generally, assume that $A$ and $B$ are graded Hopf algebras
with finite dimensional homogeneous components and that $\sigma$ is compatible with the gradings. Then the quotient Hopf
algebras $A/I_A$ and $B/I_B$ are also graded and can be identified via $\sigma$ with the duals of each other. Suppose
also that the action $\phi$ respects the grading so does the quotient $\Bar{\phi}:\pi \to \auth(A/I_A)$. In this case,
there exists a unique action $\pi \to \auth(B/I_B)$ which is $(\Bar{\sigma},\Bar{\phi})$-compatible, where
$\Bar{\sigma}\co A/I_A \times B/I_B \to \Bbbk$ is the induced Hopf pairing. Then the Hopf \p coalgebra
$D(A/I_A,B/I_B;\Bar{\sigma},\Bar{\phi})$ is quasitriangular by the same construction as in Theorem~\ref{quasidouble}.
\end{rems}
\begin{proof}
Fix basis $(e_i)$ of $A$ and $(f_i)$ of $B$ such that $\sigma(e_i,f_j)=\delta_{i,j}$ (such basis always exist since
$\sigma$ is non-degenerate). Note that $x=\sum_i \sigma(x,f_i) \, x$ and $y=\sum_i \sigma(e_i,y) \, y$ for any $x \in A$
and $y \in B$.

Recall that, since $\sum_i e_i \otimes 1_B \otimes 1_A \otimes f_i$ is the \R matrix of the usual quantum double
$D(A,B,\sigma,\id_A)$, we have
\begin{align}
& \sum_{i,j} S(e_i) e_j \otimes  f_i f_j=1_A \otimes 1_B, \label{UQD1} \\
& \sum_i e_i \otimes f_{i(1)} \otimes f_{i(2)}=\sum_{i,j} e_ie_j \otimes  f_j \otimes f_i, \label{UQD2} \\
& \sum_i e_{i(1)} \otimes e_{i(2)} \otimes f_i=\sum_{i,j} e_i \otimes e_j \otimes  f_i f_j. \label{UQD3}
\end{align}

Let $\al,\be \in \pi$. From \eqref{UQD1} and since $A$ (resp.\@ $B$) can be viewed as a subalgebra of $D(A,B;\sigma,
\phi_\al)$ (resp.\@ $D(A,B;\sigma, \phi_\be)$) via $a \mapsto a \otimes 1_B$ (resp.\@ $b \mapsto 1_A \otimes b$), we get
that $R_{\al,\be}$ is invertible in $D(A,B;\sigma, \phi_\al) \otimes D(A,B;\sigma, \phi_\be)$ with inverse
\begin{equation*}
R^{-1}_{\al,\be}= \sum_i S(e_i) \otimes 1_B \otimes 1_A \otimes f_i.
\end{equation*}

Let $a \in A$, $b \in B$ and $\al,\be \in \pi$. For all $x \in A$, we have that
\begin{eqnarray*}
\lefteqn{(\id_{A\otimes B \otimes A} \otimes \sigma(x, \cdot) )(R_{\al,\be} \cdot \cp{\al}{\be}(a \otimes b) )}&&\\
 & = & \sideset{}{_i} \sum \sigma(\phi_\be(a_{(2)}),S(f_{i(1)})) \, \sigma(a_{(4)},f_{i(3)}) \, \sigma(x,f_{i(2)}b_{(2)}) \; e_i
       \phi_\be(a_{(1)}) \otimes b_{(1)} \otimes a_{(3)} \\
 & = & \sideset{}{_i} \sum \sigma(\phi_\be S^{-1}(a_{(2)}),f_{i(1)}) \, \sigma(a_{(4)},f_{i(3)}) \, \sigma(x_{(1)},f_{i(2)})
       \sigma(x_{(2)},b_{(2)}) \; e_i \phi_\be(a_{(1)}) \otimes b_{(1)} \otimes a_{(3)} \\
 & = & \sideset{}{_i} \sum \sigma(a_{(4)} x_{(1)} \phi_\be S^{-1}(a_{(2)}),f_i) \, \sigma(x_{(2)},b_{(2)}) \;
       e_i \phi_\be(a)_{(1)} \otimes b_{(1)} \otimes a_{(3)} \\
 & = & \sigma(x_{(2)},b_{(2)}) \; a_{(4)} x_{(1)} \phi_\be(S^{-1}(a_{(2)}) a_{(1)}) \otimes b_{(1)} \otimes a_{(3)} \\
 & = & \sigma(x_{(2)},b_{(2)}) \; a_{(2)} x_{(1)} \otimes b_{(1)} \otimes a_{(1)},
\end{eqnarray*}
and, since $\displaystyle x_{(1)} \otimes x_{(2)} \otimes x_{(2)} \otimes x_{(2)}= \sum_i \sigma(x_{(2)},f_i) \; x_{(1)}
\otimes e_{i(1)} \otimes e_{i(2)} \otimes e_{i(3)}$,
\begin{eqnarray*}
\lefteqn{(\id_{A\otimes B \otimes A} \otimes \sigma(x, \cdot) )(\sigma_{\be,\al} (\varphi_\ali\otimes \id_{H_\al})
         \cp{\al \be \ali}{\al}(a \otimes b) \cdot R_{\al,\be} )}&&\\
 & = & \sideset{}{_i} \sum \sigma(\phi_\al(e_{i(1)}),S(b_{(2)})) \, \sigma(e_{i(3)},b_{(4)}) \, \sigma(x,\phi_\ali^*(b_{(1)}) f_i)
       \; a_{(2)} e_{i(2)} \otimes b_{(3)} \otimes a_{(1)} \\
 & = & \sideset{}{_i} \sum \sigma(\phi_\al(e_{i(1)}),S(b_{(2)})) \, \sigma(e_{i(3)},b_{(4)}) \, \sigma(\phi_\al(x_{(1)}),b_{(1)})
       \, \sigma(x_{(2)},f_i) \; a_{(2)} e_{i(2)} \otimes b_{(3)} \otimes a_{(1)} \\
 & = & \sigma(\phi_\al(x_{(2)}),S(b_{(2)})) \, \sigma(x_{(4)},b_{(4)}) \, \sigma(\phi_\al(x_{(1)}),b_{(1)})
       \; a_{(2)} x_{(3)} \otimes b_{(3)} \otimes a_{(1)} \\
 & = & \sigma(\phi_\al(x_{(1)}),b_{(1)}S(b_{(2)})) \, \sigma(x_{(3)},b_{(4)})
       \; a_{(2)} x_{(2)} \otimes b_{(3)} \otimes a_{(1)} \\
 & = & \sigma(x_{(2)},b_{(2)}) \; a_{(2)} x_{(1)} \otimes b_{(1)} \otimes a_{(1)}.
\end{eqnarray*}
Hence, since $\sigma$ is non-degenerate, Axiom~\eqref{QT1} is satisfied.

Let us verify Axioms~\eqref{QT4}. Let $\al,\be,\ga \in \pi$. Since $\phi^*$ is $(\sigma,\phi)$-compatible (by
definition), the basis $(\phi_\be(e_i))_i$ of $A$ and $(\phi_\be^*(f_i))_i$ of $B$ satisfy
$\sigma(\phi_\be(e_i),\phi_\be^*(e_j))=\sigma(e_i,f_j)=\delta_{i,j}$. Therefore we have that
\begin{equation*}
(\varphi_\be \otimes \varphi_\be)(R_{\al,\ga})
 = \sum_i \phi_\be(e_i) \otimes 1_B \otimes 1_A \otimes \phi^*_\be(f_j)
 = R_{\be \al \bei, \be \ga \bei}.
\end{equation*}

Finally, let us check Axioms~\eqref{QT2} and \eqref{QT3}. Let $\al,\be,\ga \in \pi$. Using \eqref{UQD2}, we have
\begin{eqnarray*}
(\id_{D(A,B;\sigma, \phi_\al)} \otimes \cp{\be}{\ga})(R_{\al,\be \ga})
 & = & \sideset{}{_i}\sum e_i \otimes 1_B \otimes 1_A \otimes f_{i(1)} \otimes 1_A \otimes f_{i(2)} \\
 & = & \sideset{}{_{i,j}} \sum e_ie_j \otimes  1_B \otimes 1_A \otimes f_j \otimes 1_A \otimes f_i\\
 & = & (R_{\al,\ga})_{1 \be 3} \cdot (R_{\al,\be})_{12 \ga}.
\end{eqnarray*}
Using \eqref{UQD3} and \eqref{QT4}, we have
\begin{eqnarray*}
\lefteqn{(\cp{\al}{\be} \otimes \id_{D(A,B;\sigma, \phi_\ga)})(R_{\al \be, \ga})}&&\\
 & = & \sideset{}{_i}\sum \phi_\be(e_{i(1)}) \otimes 1_B \otimes e_{i(2)} \otimes 1_B \otimes 1_A \otimes f_i \\
 & = & \sideset{}{_{i,j}} \sum \phi_\be(e_i) \otimes  1_B \otimes e_j \otimes 1_B \otimes 1_A \otimes f_i f_j \\
 & = & [(\varphi_\be \otimes \id_{D(A,B;\sigma, \phi_\ga)})(R_{\bei\al \be, \ga})]_{1 \be 3} \cdot (R_{\be,\ga})_{\al 23}\\
 & = & [(\id_{D(A,B;\sigma, \phi_\al)} \otimes \varphi_{\bei})(R_{\al, \be \ga \bei})]_{1 \be 3} \cdot (R_{\be,\ga})_{\al 23}.
\end{eqnarray*}
This completes the proof of the quasitriangularity of $D(A,B;\sigma, \phi)$.
\end{proof}

The following corollary is a direct consequence of Corollary~\ref{crosednondegen} and Theorem~\ref{quasidouble}.
\begin{cor}\label{Leboncor}
Let $A$ be a finite-dimensional Hopf algebra and $\phi \co \pi \to \auth (A)$ be an action of $\pi$ on $A$ by Hopf
algebras automorphisms. Recall that the duality bracket $\langle, \rangle_{A \otimes A^*}$ is a non-degenerate Hopf
pairing between $A$ and $A^{*\cop}$. Then $D(A,A^{*\cop};\langle, \rangle_{A \otimes A^*}, \phi)$ is a quasitriangular
Hopf \p coalgebra.
\end{cor}

\begin{rem}
Recall that the group of Hopf automorphisms of a finite-dimensional semisimple Hopf algebra $A$ over a field of
characteristic $0$ is finite (see \cite{Rad4}). To obtain non-trivial examples of (quasitriangular) Hopf \p coalgebras
for an infinite group $\pi$ by using the twisted double method, one has to consider non-semisimple Hopf algebras (at
least in characteristic 0).
\end{rem}

\subsection{The $h$-adic case}\label{sect-hadiccase} In this subsection, we develop the \h adic variant of Hopf
group-coalgebras. A technical argument for the need of \h adic Hopf group-coalgebras is that they are necessary for a
mathematically rigorous treatment of \R matrices for quantize enveloping algebras endowed with a group action.

Recall that if $V$ is a vector space over $\C[[h]]$, the topology on $V$ for which the sets $\{ h^nV+v | \, n \in \N \}$
are a neighborhood base of $v \in V$ is called the \emph{\h adic topology}. If $V$ and $W$ are vector spaces over
$\C[[h]]$, we shall denote by $V \hotimes W$ the completion of the tensor product space $V \otimes_{\C[[h]]} W$ in the \h
adic topology. Let $V$ be a complex vector space. Then the set $V[[h]]$ of all formal power series $f=\sum_{n=0}^\infty
v_n h^n$ with coefficients $v_n \in V$ is a vector space over $\C[[h]]$ which is complete in the \h adic topology.
Furthermore, $V[[h]] \hotimes W[[h]]=(V \otimes W)[[h]]$ for any complex vector spaces $V$ and $W$.

An \emph{\h adic algebra} is a vector space $A$ over $\C[[h]]$ which is complete in the \h adic topology and endowed with
a $\C[[h]]$-linear map $m\co A \hotimes A \to A$ and an element $1 \in A$ satisfying $m(\id_A \hotimes m)=m(m \hotimes
\id_A)$ and $m(a \hotimes 1)=a=m(1 \hotimes a)$ for all $a \in A$.

By an \emph{\h adic Hopf \p coalgebra}, we shall mean a family $H=\{H_\al\}_{\al \in \pi}$ of \h adic algebras which is
endowed with \h adic algebra homomorphisms $\cp{\al}{\be}\co H_{\al\be} \to H_\al \hotimes H_\be$ ($\al,\be \in \pi$) and
$\varepsilon \co A \to \C[[h]]$ satisfying \eqref{coass} and \eqref{counit}, and with $C[[h]]$-linear maps $S_\al \co
H_\al \to H_\ali$ ($\al \in \pi$) satisfying \eqref{antipode}. In the previous axioms, one has to replace the algebraic
tensor products $\otimes$ by the \h adic completions $\hotimes$.

The notions of crossed and quasitriangular \h adic Hopf \p coalgebras can be defined similarly as in
Sections~\ref{deficro} and \ref{quasitrig}.

The definitions of Section~\ref{sect-double} and Theorem~\ref{Hopfdouble} carry over almost verbatim to \h adic Hopf
algebras. The only modifications are that $\sigma \co A \hotimes B \to \C[[h]]$ is $\C[[h]]$-linear and that the algebra
$D(A,B; \sigma, \phi)$, where $\phi$ is an \h adic Hopf endomorphism of $A$, is built over the completion $A \hotimes B$
of $A \otimes B$ in the \h adic topology. The reasoning of the proof of Theorem~\ref{quasidouble} give the following \h
adic version.

\begin{thm}\label{doublehadicrmat}
Let $\sigma \co A \hotimes B \to \C[[h]]$ be an \h adic Hopf pairing between two \h adic Hopf algebras $A$ and $B$, and
$\phi \co \pi \to \auth (A)$ be an action of $\pi$ on $A$ by \h adic Hopf automorphisms. Suppose that $\sigma$ is
non-degenerate. Let $(e_i)_i$ and $(f_i)_i$ be dual basis of the vector spaces $A$ and $B$, respectively, with respect to
the form $\sigma$. If $R_{\al,\be}=\sum_i (e_i \otimes 1_B ) \otimes (1_A \otimes f_i)$ belongs to the \h adic completion
$D(A,B;\sigma, \phi_\al) \hotimes D(A,B;\sigma, \phi_\be)$, then $R=\{R_{\al,\be}\}_{\al,\be \in \pi}$ is a \R matrix of
the crossed \h adic Hopf \p coalgebra $D(A,B;\sigma, \phi)=\{ D(A,B;\sigma, \phi_\al)\}_{\al \in \pi}$.
\end{thm}

\section{The case of algebras of finite groups}\label{sect-finitegroup}

Let $G$ be a finite group. In this section, we describe the Hopf $G$-coalgebras obtained by the twisted double method
from the Hopf algebra $\Bbbk [G]$.

Recall that the Hopf algebra structure of the (finite-dimensional) \kt algebra $\Bbbk [G]$ of $G$ is given by
$\Delta(g)=g \otimes g$, $\varepsilon(g)=1$ and $S(g)=g^{-1}$ for all $g \in G$. The dual of $\Bbbk [G]$ is the Hopf
algebra $F(G)=\Bbbk^G$ of functions $G \to \Bbbk$. It has a basis $(e_g\co G \to \kk)_{g \in G}$ defined by $e_g(h)=
\delta_{g,h}$ where $\delta_{g,g}=1$ and $\delta_{g,h}=0$ if $g \neq h$. The structure maps of $F(G)$ are given by $e_g
e_h = \delta_{g,h} \, e_g$, $1_{F(G)}  = \sum_{g \in \pi} e_g$, $\Delta(e_g) =\sum_{xy=g} e_x \otimes e_y$,
$\varepsilon(e_g) =\delta_{g,1}$ and $S(e_g) =e_{g^{-1}}$ for any $g,h \in G$.

Set $\phi \co \Bbbk [G] \to \auth(\Bbbk [G])$ defined by $\phi_g(h)=ghg^{-1}$. It is a well defined group homomorphism
(since any $g \in G$ is grouplike in $\Bbbk[G]$). By Corollary~\ref{Leboncor}, this data leads to a quasitriangular Hopf
$G$-coalgebra $D(\Bbbk [G], F(G); \langle, \rangle_{\Bbbk [G] \times F(G)}, \phi)$ which will be denoted by
$D_G(G)=\{D_\al(G)\}_{\al \in G}$.

Let us describe $D_G(G)$ more precisely. For any $\al \in G$, the algebra structure of $D_\al(G)$, which is equal to
$\Bbbk [G] \otimes F(G)$ as a \kt space, is given by
\begin{align*}
 & (g \otimes e_h) \cdot (g' \otimes e_{h'})=\delta_{\al g' \ali, h^{-1} g' h'} \, gg' \otimes e_{h'} \text{\quad for all
   $g,g',h,h' \in G$,} \\
 & 1_{D_\al(G)}=\sum_{g \in G} 1 \otimes e_g.
\end{align*}
The structure maps of $D_G(G)$ are given, for any $\al, \be \in G$ and $g,h \in G$, by
\begin{align*}
 & \cp{\al}{\be}(g \otimes e_h)= \sum_{xy=h} \be g \bei \otimes e_y \otimes g \otimes e_x, \\
 & \varepsilon(g \otimes e_h)=\delta_{h,1},\\
 & S_\al(g \otimes e_h)=\al g^{-1} \ali \otimes e_{\al g \ali h^{-1} g^{-1}},\\
 & \varphi_\al(g \otimes e_h)=\al g \ali \otimes e_{\al h \ali}.
\end{align*}
The crossed Hopf $G$-coalgebra $D_G(G)$ is quasitriangular and furthermore ribbon with \R matrix and twist given, for any
$\al,\be \in G$, by
\begin{equation*}
R_{\al,\be}= \sum_{g,h \in G} g \otimes e_h \otimes 1 \otimes e_g \text{\quad and \quad} \theta_\al=\sum_{g \in G} \ali g
\al \otimes e_g.
\end{equation*}
Note that $\ds \theta_\al^n=\sum_{g \in G} \al^{-n} (g \al)^n \otimes e_g$ for any $n \in \mathbb{Z}$.

\section{Example of a quasitriangular Hopf $GL_n(k)$-coalgebra}\label{sect-gnlk}

In this section, $\Bbbk$ is a field whose characteristic is not $2$. Fix a positive integer $n$. We use a (finite
dimensional) Hopf algebra whose group of automorphisms is known to be the group $\mathrm{GL}_n(\Bbbk)$ of invertible
$n\times n$-matrices with coefficients in $\Bbbk$ (see~\cite{Rad4}) to derive an example of a quasitriangular Hopf
$\mathrm{GL}_n(\Bbbk)$-coalgebra.

\begin{propdef}
For $\al=(\al_{i,j}) \in \mathrm{GL}_n(\Bbbk)$, let $\mathcal{A}_n^\al$ be the $\C$-algebra generated $g$, $x_1, \dots
x_n$, $y_1, \dots, y_n$, subject to the following relations
\begin{align}
 & g^2=1, \quad x_1^2= \dots = x_n^2=0,  \quad g x_i=-x_ig, \quad  x_i x_j=-x_j x_i, \label{Andef1}\\
 & y_1^2= \dots = y_n^2=0, \quad g y_i=-y_ig, \quad  y_iy_j=-y_jy_i, \label{Andef2}\\
 & x_i y_j-y_j x_i=(\delta_{i,j} -\al_{j,i})\, g, \label{Andef3}
\end{align}
where $1 \leq i,j \leq n$. The family $\mathcal{A}_n=\{\mathcal{A}_n^\al\}_{\al \in \mathrm{GL}_n(\Bbbk)}$ has a
structure of a crossed Hopf $\mathrm{GL}_n(\Bbbk)$-coalgebra given, for any $\al=(\al_{i,j})\in\mathrm{GL}_n(\Bbbk)$,
$\be=(\be_{i,j}) \in\mathrm{GL}_n(\Bbbk)$, and $1 \leq i \leq n$, by
\begin{align}
 & \cp{\al}{\be} (g)=g \otimes g, \quad \varepsilon(g)=1, \quad S_\al(g)=g,\label{Andef4}\\
 & \cp{\al}{\be} (x_i)=1 \otimes x_i + \sum_{k=1}^n \be_{k,i} \, x_k \otimes g, \quad \varepsilon(x_i)=0,
    \quad S_\al(x_i)=\sum_{k=1}^n \al_{k,i} \, g x_k,\label{Andef5}\\
 & \cp{\al}{\be} (y_i)=y_i \otimes 1 + g \otimes y_i, \quad \varepsilon(y_i)=0, \quad S_\al(y_i)=-hy_i,\label{Andef6}\\
 & \pa(g)=g, \quad \pa(x_i)=\sum_{k=1}^n \al_{k,i} \, x_k, \quad \pa(y_i)=\sum_{k=1}^n \Tilde{\al}_{k,i}\, y_k, \label{Andef7}
\end{align}
where $(\Tilde{\al}_{i,j})=\al^{-1}$. Moreover $\mathcal{A}_n$ is quasitriangular with \R matrix given, for any $\al,\be
\in \mathrm{GL}_n(\Bbbk)$, by
\begin{equation*}
R_{\al,\be} = \frac{1}{2} \sum_{S\subseteq [n]} x_S \otimes y_S + x_S \otimes g y_S + g x_S \otimes y_S +
         g x_S \otimes g y_S.
\end{equation*}
Here $[n]=\{1, \dots n \}$, $x_\emptyset=1$, $y_\emptyset=1$, and, for a nonempty subset $S$ of $[n]$, we let
$x_S=x_{i_1} \cdots x_{i_s}$ and $y_S=y_{i_1} \cdots y_{i_s}$ where $i_1 < \cdots < i_s$ are the elements of $S$.
\end{propdef}

\begin{rem}\label{remnoteq}
From relations \eqref{Andef3}, it can be shown that the algebras $\mathcal{A}_n^\al$ and $\mathcal{A}_n^\be$ are in
general not isomorphic when $\al,\be \in \mathrm{GL}_n(\Bbbk)$ are such that $\al \neq \be$.
\end{rem}

\begin{proof}
Let $A_n$ be the \kt algebra generated by $g,x_1, \dots, x_n$ which satisfy the relations \eqref{Andef1}. The algebra
$A_n$ is $2^{n+1}$-dimensional and has a Hopf algebra structure given by
\begin{align*}
& \Delta(g) = g \otimes g, && \varepsilon(g)=1, && S(g)=g, \\
& \Delta(x_i)=x_i \otimes g + 1 \otimes x_i, && \varepsilon(x_i)=0, && S(x_i)=g x_i.
\end{align*}
Radford \cite{Rad4} showed that the group of Hopf automorphisms of $A_n$ is isomorphic to the group
$\mathrm{GL}_n(\Bbbk)$ of invertible $n\times n$-matrices with coefficients in $\Bbbk$. This group automorphism $\phi\co
\mathrm{GL}_n(\Bbbk) \to \auth(A_n)$ is given, for any $\al=(\al_{i,j}) \in \mathrm{GL}_n(\Bbbk)$, by
\begin{equation*}
\phi_\al(g)=g \text{\quad and \quad} \phi_\al(x_i)=\sum_{k=1}^n \al_{k,i} \, x_k.
\end{equation*}
The Hopf algebra $B_n=A_n^\cop$ is the \kt algebra generated by the symbols $h,y_1, \dots, y_n$ which satisfy the
relations $h^2=1$ and \eqref{Andef2} and its Hopf algebra structure is given by
\begin{align*}
& \Delta(h) = h \otimes h, && \varepsilon(h)=1, && S(h)=h, \\
& \Delta(y_i)=y_i \otimes 1 + h \otimes y_i, && \varepsilon(y_i)=0, && S(y_i)=-hy_i.
\end{align*}
Let us denote the cardinality of a set $T$ by $|T|$. The elements $g^k x_S$ (resp.\@ $h^ky_S$), where $k \in \{0,1\}$ and
$S \subseteq [n]$, form a basis for $A_n$ (resp.\@ $B_n$). Since $\Delta$ is multiplicative, it follows that
\begin{align}
& \Delta(g^kx_S)=\sum_{T \subseteq S} \lambda_{T,S} \; g^k x_T \otimes g^{k+|T|}x_{S \setminus T} \text{\quad and \quad}
\label{comultAn} \\
& \Delta(h^ky_S)=\sum_{T \subseteq S} \lambda_{T,S} \; h^{k+|T|}y_{S \setminus T} \otimes h^k y_T, \label{comultBn}
\end{align}
where $\lambda_{T,S}=\pm 1$ and $\lambda_{\emptyset,S}=1=\lambda_{S,S}$.

By Section~\ref{hopfpar}, there exists a (unique) Hopf pairing $\sigma \co A_n \times B_n \to \Bbbk$ such that, for any
$1 \leq i,j \leq n$,
\begin{equation*}
\sigma(g,h)=-1, \quad \sigma(g,y_j)=\sigma (x_i,h)=0 \text{\quad and \quad} \sigma(x_i,y_j)=\delta_{i,j}.
\end{equation*}
Using \eqref{comultAn} and \eqref{comultBn}, one gets (by induction on $|S|$) that
\begin{equation*}
\sigma(g^kx_S,h^ly_T)=(-1)^{kl} \, \delta_{S,T}
\end{equation*}
for any $k,l \in \{0,1\}$ and $S,T \subseteq [n]$, where $\delta_{S,S}=1$ and $\delta_{S,T}$ if $S \neq T$. Set
$z_0=(1+h)/2$ and $z_1=(1-h)/2$. The elements $z_k y_S$, where $k \in \{0,1\}$ and $S \subseteq [n]$, form a basis for
$B_n$ such that
\begin{equation}\label{calcsigmaAnBn}
\sigma(g^kx_S,z_ly_T)=\delta_{k,l} \, \delta_{S,T}
\end{equation}
for any $k,l \in \{0,1\}$ and $S,T \subseteq [n]$. Therefore the pairing $\sigma$ is non-degenerate. Note that this
implies that $A_n^*\cong A_n$ as a Hopf algebra.

By Theorem~\ref{quasidouble}, we get a quasitriangular Hopf $\mathrm{GL}_n(\Bbbk)$-coalgebra $D(A_n,B_n;\sigma,\phi)$.
For any $\al=(\al_{i,j}) \in \mathrm{GL}_n(\Bbbk)$, $D(A_n,B_n;\sigma,\phi_\al)$ is the algebra generated by $g$, $h$,
$x_1, \dots x_n$, $y_1, \dots, y_n$, subject to the relations $h^2=1$, \eqref{Andef1}, \eqref{Andef2} and the following
relations
\begin{align}
 & g h=hg, \quad g y_j= - y_j g, \quad h x_i = x_i h , \label{demDAn1}\\
 & x_i y_j-y_j x_i=\delta_{i,j}h -\al_{j,i}g. \label{demDAn2}
\end{align}
Indeed $D(A_n,B_n;\sigma,\phi_\al)$ is the free algebra generated by the algebras $A_n$ and $B_n$ with cross relation
\eqref{crossrel1}. Further, it suffices to require the cross relations \eqref{crossrel1} for $(1 \otimes b) \cdot (a
\otimes 1)$ with $a=g, x_i$ and $b=h,y_j$. To simplify the notations, we identify of $a$ with $a \otimes 1$ and of $b$
with $1 \otimes b$ (recall that these natural maps $A_n \hookrightarrow D(A_n,B_n;\sigma,\phi_\al)$ and $B_n
\hookrightarrow D(A_n,B_n;\sigma,\phi_\al)$ are algebra monomorphisms). For example, let $a=x_i$ and $b=y_j$. Since
$\sigma(x_i,1)=\sigma(g,y_j)=\sigma(x_i,h)=\sigma(1,y_j)=0$, relation \eqref{crossrel1} gives
\begin{equation*}
y_j x_i=\sigma(\phi_\al(x_i),y_jh)\sigma(g,1)\, g \cdot 1 + \sigma(1,h) \sigma(g,1) x_i \cdot y_j + \sigma(1,h)
\sigma(x_i,y_j) 1 \cdot h.
\end{equation*}
Inserting the values $\sigma(g,1)=\sigma(1,h)=1$, $\sigma(x_i,y_j)=\delta_{i,j}$, and
$\sigma(\phi_\al(x_i),y_jh)=-\al_{j,i}$, we get \eqref{demDAn2}.

From Theorem~\ref{Hopfdouble}, we obtain that the comultiplication $\cp{\al}{\be}$, the counit $\varepsilon$, the
antipode $S_\al$, and the crossing $\varphi_\al$ of $D(A_n,B_n;\sigma,\phi_\al)$ are given by
\begin{align}
 & \cp{\al}{\be} (g)=g \otimes g, \quad \cp{\al}{\be} (h)=h \otimes h, \label{demDAn3}\\
 & \cp{\al}{\be} (x_i)=1 \otimes x_i + \sum_{k=1}^n \be_{k,i} \, x_k \otimes g, \quad
   \cp{\al}{\be} (y_i)=y_i \otimes 1 + h \otimes y_i, \label{demDAn4}\\
 & \varepsilon(g)=\varepsilon(h)=1, \quad \varepsilon(x_i)=\varepsilon(y_i)=0, \quad S_\al(g)=g, \label{demDAn5}\\
 &  S_\al(h)=h, \quad S_\al(x_i)=\sum_{k=1}^n \al_{k,i} \, g x_k,
   \quad S_\al(y_i)=-hy_i,\label{demDAn6}\\
 & \pa(g)=g, \quad \pa(h)=h, \quad \pa(x_i)=\sum_{k=1}^n \al_{k,i} \, x_k, \quad \pa(y_i)=\sum_{k=1}^n \Tilde{\al}_{k,i}
 \, y_k,\label{demDAn7}
\end{align}
where $(\Tilde{\al}_{i,j})=\al^{-1}$.

For any $\al \in \mathrm{GL}_n(\Bbbk)$, let $I_\al$ be the ideal of $D(A_n,B_n;\sigma,\phi_\al)$ generated by $g-h$.
Using the above description of the structure maps of $D(A_n,B_n;\sigma,\phi)$, we get that $I=\{I_\al\}_{\al \in \pi}$ is
a crossed Hopf $\mathrm{GL}_n(\Bbbk)$-coideal of $D(A_n,B_n;\sigma,\phi)$. The quotient
$D(A_n,B_n;\sigma,\phi)/I=\{D(A_n,B_n;\sigma,\phi_\al)/I_\al \}_{\al \in \mathrm{GL}_n(\Bbbk)}$ is is precisely
$\mathcal{A}_n=\{\mathcal{A}_n^\al\}_{\al \in \mathrm{GL}_n(\Bbbk)}$ and so the latter has a quasitriangular Hopf
$\mathrm{GL}_n(\Bbbk)$-coalgebra structure which can be described by replacing $h$ with $g$ in
\eqref{demDAn3}-\eqref{demDAn7}.

Finally, the \R matrix of $\mathcal{A}_n$ is obtained as the image under the projection maps $D(A_n,B_n;\sigma,\phi_\al)
\overset{p_\al}{\to} D(A_n,B_n;\sigma,\phi_\al)/I_\al=\mathcal{A}_n^\al$ of the \R matrix of $D(A_n,B_n;\sigma,\phi)$,
that is, using \eqref{calcsigmaAnBn},
\begin{eqnarray*}
R_{\al,\be} & = & \sum_{S\subseteq [n]} p_\al(x_S) \otimes p_\be(z_0 y_S) + p_\al(gx_S) \otimes p_\be(z_1 y_S) \\
   & = & \sum_{S\subseteq [n]} x_S \otimes (\frac{1+g}{2}) y_S + gx_S \otimes (\frac{1-g}{2}) y_S \\
   & = & \frac{1}{2} \sum_{S\subseteq [n]} x_S \otimes y_S + x_S \otimes g y_S + g x_S \otimes y_S +
         g x_S \otimes g y_S.
\end{eqnarray*}
\end{proof}

\section{Graded quantum groups}\label{sect-uq}

Let $\g$ be a finite-dimensional complex simple Lie algebra of rank $l$ with Cartan matrix $(a_{i,j})$. We let $d_i$ be
the coprime integers such that the matrix $(d_i a_{i,j})$ is symmetric. Let $q$ be a fixed non-zero complex number and
let $q_i=q^{d_i}$. Suppose that $q_i^2 \neq 1$ for $i=1,2, \dots, l$.

\begin{propdef}\label{uqgpi}
Set $\pi=(\C^*)^l$. For $\al=(\al_1, \dots,\al_l) \in \pi$, let $U^\al_q(\g)$ be the $\C$-algebra generated by $K_i^{\pm
1}$, $E_i$, $F_i$, $1 \leq i \leq l$, subject to the following defining relations:
\begin{align}
 & K_iK_j=K_jK_i, \quad K_i K_i^{-1}= K_i^{-1} K_i =1, \label{uq1}\\
 & K_i E_j = q_i^{a_{i,j}} E_j K_i, \label{uq2}\\
 & K_i F_j = q_i^{-a_{i,j}} F_j K_i, \label{uq3} \\
 & E_i F_j-F_j E_i=\delta_{i,j}\frac{\al_i K_i - K_i^{-1}}{q_i-q_i^{-1}} \label{uq4}\\
 & \sum_{r=0}^{1-a_{i,j}} (-1)^r \begin{bmatrix}  1-a_{i,j} \\ r \end{bmatrix}_{q_i} E_i^{1-a_{i,j}-r}E_j E_i^r=0 \text{\quad
   if \quad } i \neq j,\label{uq5} \\
 & \sum_{r=0}^{1-a_{i,j}} (-1)^r \begin{bmatrix}  1-a_{i,j} \\ r \end{bmatrix}_{q_i} F_i^{1-a_{i,j}-r}F_j F_i^r=0 \text{\quad
   if \quad } i \neq j. \label{uq6}
\end{align}
The family $U^\pi_q(\g)=\{U^\al_q(\g)\}_{\al \in \pi}$ has a structure of a crossed Hopf \p coalgebra given, for
$\al=(\al_1, \dots,\al_l)\in\pi$, $\be=(\be_1, \dots,\be_l) \in \pi$ and $1 \leq i \leq l$, by
\begin{align*}
 & \cp{\al}{\be} (K_i)=K_i \otimes K_i,\\
 & \cp{\al}{\be} (E_i)=\be_i E_i \otimes K_i + 1 \otimes E_i, \\
 &  \cp{\al}{\be} (F_i)=F_i \otimes 1 + K_i^{-1} \otimes F_i,\\
 & \varepsilon(K_i)=1, \quad \varepsilon(E_i)=\varepsilon(F_i)=0,\\
 & S_\al(K_i)=K_i^{-1}, \quad S_\al(E_i)=-\al_i E_i K_i^{-1}, \quad S_\al(F_i)=-K_i F_i,\\
 & \pa(K_i)=K_i, \quad \pa(E_i)=\al_i E_i, \quad \pa(F_i)=\al_i^{-1} F_i.
\end{align*}
\end{propdef}
\begin{rem}
Note that $(U^1_q(\g), \cp{1}{1}, \varepsilon, S_1)$ is the usual quantum group $U_q(\g)$.
\end{rem}

\begin{proof}
Let $U_+$ be the $\C$-algebra generated by $E_i$, $K_i^{\pm 1}$, $1 \leq i \leq l$, subject to the relations \eqref{uq1},
\eqref{uq2} and \eqref{uq5}. Let $U_-$ be the $\C$-algebra generated by $F_i$, ${K'_i}^{\pm 1}$, $1 \leq i \leq l$,
subject to the relations \eqref{uq1}, \eqref{uq3} and \eqref{uq6} where one has to replace $K_i$ with $K'_i$. The
algebras $U_+$ and $U_-$ have a Hopf algebra structure given by
\begin{align*}
 & \Delta(K_i)=K_i \otimes K_i, \quad \Delta(E_i)=E_i \otimes K_i + 1 \otimes E_i,\\
 & \varepsilon(K_i)=1, \quad \varepsilon(E_i)=0,\quad
 S(K_i)=K_i^{-1}, \quad S(E_i)=-E_i K_i^{-1},\\
 & \Delta(K'_i)=K'_i \otimes K'_i, \quad \Delta(F_i)=F_i \otimes 1 + {K'_i}^{-1} \otimes F_i,\\
 & \varepsilon(K'_i)=1, \quad \varepsilon(F_i)=0,\quad
  S(K_i)=K_i^{-1}, \quad S(F_i)=-K'_i F_i.
\end{align*}
Using the method described in Section~\ref{hopfpar}, it can be verified that there exists a (unique) Hopf pairing $\sigma
\co U_+ \times U_- \to \C$ such that
\begin{equation*}
 \sigma(E_i,F_j)=\frac{\delta_{i,j}}{q_i-q_i^{-1}}, \quad
 \sigma(E_i,K'_j)=\sigma(K_i,F_j)=0, \quad
 \sigma(K_i,K'_j)=q_i^{a_{i,j}}=q_j^{a_{j,i}}.
\end{equation*}
Let $\phi: \pi \to \auth (U_+)$ and $\psi: \pi \to \auth (U_-)$ defined, for $\be=(\be_1, \dots,\be_l) \in \pi$ and $1
\leq i \leq l$, by
\begin{equation*}
 \phi_\be(K_i)=K_i, \quad \phi_\be(E_i)= \be_i \, E_i, \quad \psi_\be(K'_i)=K'_i, \quad \psi_\be(F_i)= \be_i^{-1} \, F_i.
\end{equation*}
It is straightforward to verify that $\psi$ is $(\sigma,\phi)$-compatible. By Lemma~\ref{crosseddegen}, we can consider
the crossed Hopf \p coalgebra $D(U_+,U_-;\sigma,\phi)=\{D(U_+,U_-;\sigma,\phi_\al)\}_{\al \in \pi}$.

Now, for any $\al \in \pi$, $D(U_+,U_-;\sigma,\phi_\al)$ is the algebra generated by $K_i^{\pm 1}$, ${K'_i}^{\pm 1}$,
$E_i$, $F_i$,  where $1 \leq i \leq l$, subject to the relations \eqref{uq1}, \eqref{uq2}, \eqref{uq5}, the relations
\eqref{uq1}, \eqref{uq3}, \eqref{uq6} where one has to replace $K_i$ with $K'_i$, and the following relations
\begin{align}
 & K_i K'_j=K'_j K_i, \quad K_i F_j=q_i^{-a_{i,j}} F_j K_i, \quad K'_i E_j = q_i^{a_{i,j}} E_j K'_i, \label{uqp1}\\
 & E_i F_j-F_j E_i=\delta_{i,j}\frac{\al_i K_i - {K'_i}^{-1}}{q_i-q_i^{-1}}. \label{uqp2}
\end{align}
Indeed $D(U_+,U_-;\sigma,\phi_\al)$ is the free algebra generated by the algebras $U_+$ and $U_-$ with cross relation
\eqref{crossrel1}. Further, it suffices to require the cross relations \eqref{crossrel1} for $(1 \otimes b) \cdot (a
\otimes 1)$ with $a=K_i, E_i$ and $b=K'_i,F_i$. To simplify the notations, we identify of $a$ with $a \otimes 1$ and of
$b$ with $1 \otimes b$ (recall that these natural maps $U_+ \hookrightarrow D(U_+,U_-;\sigma,\phi_\al)$ and $U_-
\hookrightarrow D(U_+,U_-;\sigma,\phi_\al)$ are algebra monomorphisms). For example, let $a=E_i$ and $b=F_j$. Since
$\sigma(E_i,1)=\sigma(K_i,F_j)=\sigma(E_i,{K'_j}^{-1})=\sigma(1,F_j)=0$, relation \eqref{crossrel1} gives
\begin{equation*}
F_j E_i=\sigma(\al_i E_i,S(F_j)) \sigma(K_i,1)\, K_i + \sigma(1,K'_j) \sigma(K_i,1) E_iF_j + \sigma(1,K'_j)
\sigma(E_i,F_j) \,{K'_j}^{-1}.
\end{equation*}
Inserting the values $\sigma(K_i,1)=\sigma(1,K'_j)=1$, $\sigma(E_i,F_j)=\delta_{i,j}(q_i-q_i^{-1})^{-1}$ and
$\sigma(E_i,S(F_j))=-\delta_{i,j}(q_i-q_i^{-1})^{-1}$, we get \eqref{uqp2}.

From Theorem~\ref{Hopfdouble}, we obtain that the comultiplication $\cp{\al}{\be}$, the counit $\varepsilon$, the
antipode $S_\al$, and the crossing $\varphi_\al$ of $D(U_+,U_-;\sigma,\phi)$ are given, for $1 \leq i \leq l$, by
\begin{align}
 & \cp{\al}{\be} (K_i)=K_i \otimes K_i, \quad \cp{\al}{\be} (K'_i)=K'_i \otimes K'_i, \label{uqpmap1}\\
 & \cp{\al}{\be} (E_i)=\be_i E_i \otimes K_i + 1 \otimes E_i, \quad
   \cp{\al}{\be} (F_i)=F_i \otimes 1 + {K'_i}^{-1} \otimes F_i, \label{uqpmap2}\\
 & \varepsilon(K_i)=\varepsilon(K'_i)=1, \quad \varepsilon(E_i)=\varepsilon(F_i)=0, \quad S_\al(K_i)=K_i^{-1}, \label{uqpmap3}\\
 &  S_\al(K'_i)={K'_i}^{-1}, \quad S_\al(E_i)=-\al_i E_i K_i^{-1},
   \quad S_\al(F_i)=-K'_i F_i,\label{uqpmap4}\\
 & \pa(K_i)=K_i, \quad \pa(K'_i)=K'_i, \quad \pa(E_i)=\al_i E_i, \quad \pa(F_i)=\al_i^{-1} F_i.\label{uqpmap5}
\end{align}

Finally, for any $\al \in \pi$, let $I_\al$ be the ideal of $D(U_+,U_-;\sigma,\phi_\al)$ generated by $K_i-K'_i$ and
$K_i^{-1}-{K'_i}^{-1}$, where $1\leq i \leq l$. Using the above description of the structure maps of
$D(U_+,U_-;\sigma,\phi)$, we get that $I=\{I_\al\}_{\al \in \pi}$ is a crossed Hopf \p coideal of
$D(U_+,U_-;\sigma,\phi)$. The quotient $D(U_+,U_-;\sigma,\phi)/I=\{D(U_+,U_-;\sigma,\phi_\al)/I_\al \}_{\al \in \pi}$ is
precisely $U_q^\pi(\g)=\{U_q^\al\}_{\al \in \pi}$ and so the latter has a crossed Hopf \p coalgebra structure which can
be described by replacing $K'_i$ with $K_i$ in \eqref{uqpmap1}-\eqref{uqpmap5}.
\end{proof}

\section{$h$-adic graded quantum groups}\label{sect-uh}

Let $\g$ be a finite-dimensional complex simple Lie algebra of rank $l$ with Cartan matrix $(a_{i,j})$. We let $d_i$ be
the coprime integers such that the matrix $(d_i a_{i,j})$ is symmetric.

\begin{propdef}\label{uqhpi}
Set $\pi=\C[[h]]^l$. For $\al=(\al_1, \dots,\al_l) \in \pi$, let $U^\al_h(\g)$ be the algebra over $\C[[h]]$
topologically generated by the elements $H_i$, $E_i$, $F_i$, $1 \leq i \leq l$, subject to the following defining
relations:
\begin{align}
 & [H_i,H_j]=0, \label{uh1}\\
 & [H_i,E_j]= a_{ij} E_j, \label{uh2}\\
 & [H_i,F_j]= -a_{ij} F_j, \label{uh3} \\
 & [E_i,F_j]=\delta_{i,j}\frac{e^{d_i h\al_i} e^{d_i hH_i} - e^{-d_i hH_i}}{e^{d_i h}-e^{-d_i h}} \label{uh4}\\
 & \sum_{r=0}^{1-a_{i,j}} (-1)^r \begin{bmatrix}  1-a_{i,j} \\ r \end{bmatrix}_{e^{d_i h}} E_i^{1-a_{i,j}-r}E_j E_i^r=0 \quad
   (i \neq j),\label{uh5} \\
 & \sum_{r=0}^{1-a_{i,j}} (-1)^r \begin{bmatrix}  1-a_{i,j} \\ r \end{bmatrix}_{e^{d_i h}} F_i^{1-a_{i,j}-r}F_j F_i^r=0 \quad
   (i \neq j). \label{uh6}
\end{align}
The family $U^\pi_q(\g)=\{U^\al_q(\g)\}_{\al \in \pi}$ has a structure of a crossed Hopf \p coalgebra given, for
$\al=(\al_1, \dots,\al_l)\in\pi$, $\be=(\be_1, \dots,\be_l) \in \pi$ and $1 \leq i \leq l$, by
\begin{align*}
 & \cp{\al}{\be} (H_i)=H_i \otimes 1 + 1 \otimes H_i,\\
 & \cp{\al}{\be} (E_i)=e^{d_i h \be_i} E_i \otimes e^{d_i h H_i} + 1 \otimes E_i, \\
 &  \cp{\al}{\be} (F_i)=F_i \otimes 1 + e^{-d_i h H_i} \otimes F_i,\\
 & \varepsilon(H_i)=\varepsilon(E_i)=\varepsilon(F_i)=0,\\
 & S_\al(H_i)=-H_i, \quad S_\al(E_i)=-e^{d_i h \al_i} E_i e^{-d_i h H_i},
   \quad S_\al(F_i)=-e^{d_i h H_i} F_i,\\
 & \pa(H_i)=H_i, \quad \pa(E_i)=e^{d_i h \al_i} E_i, \quad \pa(F_i)=e^{-d_i h \al_i} F_i.
\end{align*}
\end{propdef}
\begin{rems}
{(a)} $(U^0_h(\g), \cp{0}{0}, \varepsilon, S_0)$ is the usual quantum group $U_h(\g)$.\\
\indent {(b)} The element $e^{d_ih}-e^{-d_ih}\in \C[[h]]$ is not invertible in $\C[[h]]$, because the constant term is
zero. But the expression oh the right hand side of \eqref{uh4} is a formal power series $\sigma_n p_n(H_i) h^n$ with
certain polynomials $p_n(H_i)$, and so it is a well-defined element of the $h$-adic algebra generated by $E_i$, $F_i$,
$H_i$.
\end{rems}

\begin{proof}
Let $U_+$ be the $h$-adic algebra generated by $H_i$, $E_i$, $1 \leq i \leq l$, subject to the relations \eqref{uh1},
\eqref{uh2} and \eqref{uh5}. Let $U_-$ be the $h$-adic algebra generated by $H'_i$, $F_i$, $1 \leq i \leq l$, subject to
the relations \eqref{uh1}, \eqref{uh3} and \eqref{uh6} where one has to replace $H_i$ with $H'_i$. The algebras $U_+$ and
$U_-$ have a $h$-adic Hopf algebra structure given by
\begin{align*}
 & \Delta(H_i)=H_i \otimes 1 + 1 \otimes H_i, \quad \Delta(E_i)=E_i \otimes e^{d_ihH_i} + 1 \otimes E_i,\\
 & \varepsilon(H_i)=\varepsilon(E_i)=0,\quad
 S(H_i)=-H_i, \quad S(E_i)=-E_i e^{-d_ihH_i},\\
 & \Delta(H'_i)=H'_i \otimes 1 + 1 \otimes H'_i, \quad \Delta(F_i)=F_i \otimes 1 + e^{-d_ihH'_i} \otimes F_i,\\
 & \varepsilon(H'_i)=\varepsilon(F_i)=0,\quad
  S(H'_i)=-H'_i, \quad S(F_i)=-e^{d_ihH'_i} F_i.
\end{align*}

Let us consider the $h$-adic Hopf algebra $\widetilde{U}_-=\C[[h]]1+hU_-$. The elements $\widetilde{H}'_i=hH'_i$ and
$\widetilde{F}_i=h F_i$ belong to $\widetilde{U}_-$ and satisfy
\begin{equation*}
 [\widetilde{H}_i,\widetilde{F}_j]= -ha_{ij} \widetilde{F}_j, \quad
 \Delta(\widetilde{H}'_i)=\widetilde{H}'_i \otimes 1 + 1 \otimes \widetilde{H}'_i,
 \quad \Delta(\widetilde{F}_i)=\widetilde{F}_i \otimes 1 + e^{-d_i\widetilde{H}'_i} \otimes \widetilde{F}_i.
\end{equation*}
The element $e^{-d_i\widetilde{H}'_i}=1 + \sum_{k\geq 1} \frac{1}{k!} (-d_ih)^kH_i^k$ is also in $\widetilde{U}_-$. Note
that $e^{-d_i\widetilde{H}'_i}$ is not in the $h$-adic subalgebra of $\widetilde{U}_-$ generated by $\widetilde{H}'_i$.

Using the method described in Section~\ref{hopfpar}, it can be verified that there exists a (unique) Hopf pairing $\sigma
\co U_+ \times \widetilde{U}_- \to \C[[h]]$ such that
\begin{equation*}
 \sigma(H_i,\widetilde{H}'_j)=d_i^{-1}a_{j,i}, \quad
 \sigma(H_i,\widetilde{F}_j)=\sigma(E_i,\widetilde{H}'_j)=0, \quad
 \sigma(E_i,\widetilde{F}_j)=\frac{\delta_{i,j}\,h}{e^{d_ih}- e^{-d_ih}}.
\end{equation*}
Let $\phi: \pi \to \auth (U_+)$ and $\psi: \pi \to \auth (\widetilde{U}_-)$ defined, for $\al=(\al_1, \dots,\al_l) \in
\pi$ and $1 \leq i \leq l$, by
\begin{equation*}
 \phi_\al(H_i)=H_i, \quad \phi_\al(E_i)= e^{d_ih \al_i} \, E_i, \quad \psi_\al(\widetilde{H}'_i)=\widetilde{H}'_i,
 \quad \psi_\be(\widetilde{F}_i)= e^{-d_ih \al_i} \, \widetilde{F}_i.
\end{equation*}
It is straightforward to verify that $\psi$ is $(\sigma,\phi)$-compatible. By the $h$-adic version of
Lemma~\ref{crosseddegen}, we can consider the crossed $h$-adic Hopf \p coalgebra
$D(U_+,\widetilde{U}_-;\sigma,\phi)=\{D(U_+,\widetilde{U}_-;\sigma,\phi_\al)\}_{\al \in \pi}$ whose structure can be
explicitly described as in the proof of Proposition~\ref{uqgpi}.

For any $\al \in \pi$, let $I_\al$ be the $h$-adic ideal of $D(U_+,\widetilde{U}_-;\sigma,\phi_\al)$ generated by
$\widetilde{H}'_i-hH_i$ where $1\leq i \leq l$. Using the description of the structure maps of
$D(U_+,\widetilde{U}_-;\sigma,\phi_\al)$, we get that $I=\{I_\al\}_{\al \in \pi}$ is a crossed $h$-adic Hopf \p coideal
of $D(U_+,U_-;\sigma,\phi)$. The quotient
$D(U_+,\widetilde{U}_-;\sigma,\phi)/I=\{D(U_+,\widetilde{U}_-;\sigma,\phi_\al)/I_\al \}_{\al \in \pi}$ is precisely
$U_h^\pi(\g)=\{U_h^\al\}_{\al \in \pi}$ and so the latter has a crossed $h$-adic Hopf \p coalgebra structure.
\end{proof}

It is well-know (see, e.g., \cite{KS1}) that the Hopf pairing $\sigma \co U_+ \times \widetilde{U}_- \to \C[[h]]$ is
non-degenerate and that, if $(e_i)_i$ and $(f_i)_i$ are dual basis of the vector spaces $U_+$ and $\widetilde{U}_-$ with
respect to the form $\sigma$, then $\sum_i (e_i \otimes 1) \otimes (1\otimes f_i)$ belongs to the \h adic completion
$D(U_+ ,\widetilde{U}_-;\sigma, \phi_\al) \hotimes D(U_+,\widetilde{U}_-;\sigma, \phi_\be)$. Therefore, by
Theorem~\ref{doublehadicrmat}, the crossed $h$-adic Hopf \p coalgebra $D(U_+,\widetilde{U}_-;\sigma,\phi)$ is
quasitriangular. Hence, as a quotient of $D(U_+ ,\widetilde{U}_-;\sigma, \phi)$,  $U_h^\pi(\g)$ is also quasitriangular.

For example, when $\g=\mathfrak{sl}_2$ and so $\pi=\C[[h]]$, we have that the \R matrix of
$U_h^{\C[[h]]}(\mathfrak{sl}_2)$ is given, for any $\al,\be \in \C[[h]]$, by
\begin{equation*}
R_{\al,\be}=e^{h (H \otimes H)/2} \sum_{n=0}^\infty R_n(h)\, E^n \otimes F^n,
\end{equation*}
where $R_n(h)=q^{n(n+1)/2} \frac{(1-q^{-2})^n}{[n]_q!}$ and $q=e^h$.

For $n\geq 1$, there exits a representation $\rho_n^\al:U_h^\al(\mathfrak{sl}_2) \to \mathrm{GL}(V^\al_n)$, where
$V^\al_n=\C^n$ as a vector space, given on the standard basis $(e_i)_{1\leq i \leq n}$ of $\C^n$ by
\begin{align*}
& \rho_n^\al(H)e_i=(n-2i+1-\frac{\al}{2}) \, e_i,\\
& \rho_n^\al(E)e_i=\begin{cases}e^\frac{h\al}{2}\, [n-i+1]_q \,e_{i-1} & \text{if $i>1$} \\ 0 & \text{if $i=0$} \end{cases}, \\
& \rho_n^\al(F)e_i=\begin{cases}[i]_q \,e_{i+1} & \text{if $i<n$} \\ 0 & \text{if $i=n$} \end{cases}.
\end{align*}
Together with the quasitriangularity of $U_h^{\C[[h]]}(\mathfrak{sl}_2)$, this data leads in particular to a solution of
the $\C[[h]]$-colored Yang-Baxter equation.

\bibliographystyle{amsalpha}
\bibliography{GQG}
\end{document}